\documentclass[12pt]{article}
\usepackage{amssymb}
\makeatletter
\date{19/02/2014}

\newtheorem{theorem}{Theorem}[section]

\newtheorem{definition}[theorem]{Definition}
\newtheorem{lemma}[theorem]{Lemma}

\newtheorem{corollary}[theorem]{Corollary}

\newcommand{\black}{{\blacksquare}}

\begin{document}

\title{ALMOST OVERCOMPLETE AND ALMOST OVERTOTAL SEQUENCES IN BANACH SPACES}
\author{Vladimir P. Fonf and Clemente Zanco}
\date{}
\maketitle

\footnotetext{Research of the first author was supported in part
by Israel Science Foundation, Grant \# 209/09 and by the Gruppo
Nazionale per l'Analisi Matematica, la Probabilit\`a e le loro
Applicazioni (GNAMPA) of the Istituto Nazionale di Alta Matematica
(INdAM) of Italy.} \footnotetext{Research of the second author was
supported in part by the Gruppo Nazionale per l'Analisi
Matematica, la Probabilit\`a e le loro Applicazioni (GNAMPA) of
the Istituto Nazionale di Alta Matematica (INdAM) of Italy and in
part by the Center for Advanced Studies in Mathematics at the
Ben-Gurion University of the Negev, Beer-Sheva, Israel.}

\bigskip
\bigskip
\bigskip
\bigskip
\bigskip

\begin{abstract}
The new concepts are introduced of almost overcomplete sequence in
a Banach space and almost overtotal sequence in a dual space. We
prove that any of such sequences is relatively norm-compact and we
obtain several applications of this fact.
\end{abstract}

\bigskip
\bigskip

\bigskip

{\it 2000 Mathematics Subject Classification}: Primary 46B20,
46B50; Secondary 46B45

\bigskip

 {\it Key words and phrases}: overcomplete sequence, overtotal sequence.

\vfill\eject

\maketitle

\section{Introduction}

Recall that a sequence in a Banach space $X$ is said overcomplete
in $X$ whenever the linear span of any its subsequence is dense in
$X$. It is a well-known fact that overcomplete sequences exist in
any separable Banach space. In the spirit of this notion, we
introduce the new notion of {\it overtotal sequence} and weaken
both these notions to that ones of {\it almost overcomplete
sequence} and {\it almost overtotal sequence}.

The main goal of this paper is to prove that any bounded almost
overcomplete sequence as well as any bounded almost overtotal
sequence is relatively norm-compact (section 2). We feel that
these facts provide useful tools for attacking many questions: in
section 3 several applications are presented to support this
feeling.

\medskip

Throughout the paper we use standard Geometry of Banach Spaces
terminology and notation as in \cite{JL}. In particular, $[S]$
stands for the closure of the linear span of the set $S$ and by
``subspace'' we always mean ``closed subspace''.

\medskip

Let us start by giving our three new definitions.
\begin{definition}
\label{ots} Let $X$ be a Banach space. A sequence in the dual
space $X^*$ is said to be overtotal on $X$ whenever any its
subsequence is total over $X$.
\end{definition}
If $X$ admits a total sequence $\{x_n^*\}\subset X^*$, then there
is an overtotal sequence on $X$. Indeed, put $Y=[\{ x_n^*\}]$: $Y$
is a separable Banach space, so it has an overcomplete sequence
$\{y_n^*\}.$ It is easy to see that $\{y_n^*\}$ is overtotal on
$X$.

As an easy example of an overtotal sequence, consider $X = A(D)$,
where $A(D)$ is the usual Banach disk algebra whose elements are
the holomorphic functions on the open unit disk $D$ of the plane
that admit continuous extension to $\partial D$, and $\{ x_n^*\} =
\{ z_n|_{A(D)} \}$ where $\{ z_n\}$ is any sequence of points of
$D$ converging inside $D$.

\begin{definition}
\label{ots1} A sequence in a Banach space $X$ is said to be almost
overcomplete whenever the closed linear span of any its
subsequence has finite codimension in $X$.
\end{definition}

\begin{definition}
\label{ots2} Let $X$ be a Banach space. A sequence in the dual
space $X^*$ is said to be almost overtotal on $X$ whenever the
annihilator (in $X$) of any its subsequence has finite dimension.
\end{definition}

Clearly, any overcomplete $<$ overtotal $>$ sequence is almost
overcomplete $<$ almost overtotal $>$ and the converse is not
true. It is easy to see that, if $\{ (x_n,x^*_n)\}$ is a countable
biorthogonal system, then neither $\{ x_n\}$ can be almost
overcomplete in $[\{ x_n\}]$, nor $\{ x^*_n\}$ can be almost
overtotal on $[\{ x_n\}]$. In particular, any almost overcomplete
sequence has no basic subsequence.

\section{Main results}

\bigskip

\begin{theorem}
\label{compocs} Each almost overcomplete bounded sequence in a
Banach space is relatively norm-compact.\end{theorem}

\noindent{\bf Proof.} Let $\{x_n\} $ be an almost overcomplete
bounded sequence in a (separable) Banach space $(X,|| \cdot ||)$.
Without loss of generality we may assume, possibly passing to an
equivalent norm, that the norm $|| \cdot ||$ is locally uniformly
rotund (LUR) and that $\{x_n\} $ is normalized under that norm.

First note that $\{x_n\} $ is relatively weakly compact:
otherwise, it is known (see for instance Theorem 1.3 ({\it i}) in
\cite{Si2}) that it should admit some subsequence that is a basic
sequence, a contradiction. Hence, by the Eberlein-\v{S}mulyan theorem,
$\{x_n\} $ admits some subsequence $\{x_{n_k} \}$ that weakly
converges to some point $x_0 \in B_X$. Two cases must now be
considered.

1) $||x_0|| < 1$. From $||x_{n_k} - x_0|| \geq 1 - ||x_0|| > 0$,
according to a well known result, it follows that some subsequence
$\{x_{n_{k_i}} - x_0\}$ is a basic sequence: hence ${\rm
codim}[\{x_{n_{k_{2i}}} - x_0\}] = {\rm codim}[\{x_{n_{k_{2i}}}\},
x_0] = {\rm codim}[\{x_{n_{k_{2i}}}\}] = \infty$, a contradiction.

2) $||x_0|| = 1$. Since we are working with a LUR norm, the
subsequence $\{x_{n_k} \}$ actually converges to $x_0$ in the norm
too and we are done. \ $\black$

\bigskip

\noindent {\bf Remark.} \ For overcomplete bounded sequences in
reflexive spaces this theorem has been already proved in
\cite{CFP}.

\medskip

As a first immediate consequence we get the following Corollary.

\begin{corollary}\label{corAOS}
Let $X$ be a Banach space and $\{ x_n\} \subset B_X$ be a sequence
that is not relatively norm-compact. Then there exists an
infinite-dimensional subspace $Y$ of $X^*$ such that $|\{x_n\}\cap
Y^{\top}|=\infty .$ For instance this is true for any
$\delta$-separated sequence $\{x_n\}\subset B_{X} \ (\delta
>0)$.
\end{corollary}

\begin{theorem}
\label{compots} Le $X$ be a separable Banach space. Any bounded
sequence that is almost overtotal on $X$ is relatively
norm-compact.
\end{theorem}

\noindent{\bf Proof.} Let $\{f_n \}_{n=1}^\infty \subset X^*$ be a
bounded sequence almost overtotal on $X$. Without loss of
generality, like in the proof of Theorem \ref{compocs}, we may
assume $\{f_n \} \subset S_{X^*}$. Let $\{ f_{n_k}\}$ be any
subsequence of $\{f_n \}$: since $X$ is separable, without loss of
generality we may assume that $\{ f_{n_k}\}$ weakly converges, say
to $f_0$.

Let $Z$ be a separable subspace of $X^*$ that is 1-norming for
$X$. Put $Y = [\{f_n \}_{n=0}^\infty,Z]$. Clearly $X$
isometrically embeds into $Y^*$ (we isometrically embed $X$ into
$X^{**}$ in the usual way) and $X$ is 1-norming for $Y$. By Lemma
16.3 in \cite{Si1} there is an equivalent norm $||| \cdot |||$ on
$Y$ such that, for any sequence $\{h_k \}$ and $h_0$ in $Y$,

\begin{equation}\label{leq}
h_k(x) \to h_0(x) \ \ \forall x \in X \ \ \ \  {\rm implies} \ \ \
\ |||h_0|||  \leq {\rm lim \, inf} |||h_k|||
\end{equation}

and, in addition,

\begin{equation}\label{=}
|||h_k||| \to |||h_0||| \ \ \ \ {\rm implies} \ \ \ \ |||h_k -
h_0||| \to 0.
\end{equation}

Take such an equivalent norm on $Y$ and put $h_k = f_{n_k}$ and
$h_0 = f_0$. By (\ref{=}), we are done if we prove that $|||h_k|||
\to |||h_0|||$. Suppose to the contrary that

\begin{equation}\label{not to}
|||f_{n_k}||| \not\to |||f_0|||.
\end{equation}

From (\ref{leq}) it follows that there are $\{n_{k_i} \}$ and
$\delta >0$ such that $|||f_{n_{k_i}}||| - |||f_0||| > \delta$,
that forces $||| f_{n_{k_i}} - f_0||| > \delta$ for $i$ big
enough. By \cite{JR}, Theorem III.1, it follows that some
subsequence $\{ f_{n_{k_{i_m}}} - f_0 \}_{m=1}^\infty$ is a
$w^*$-basic sequence (remember that $Y \subset X^*$, $X$ is
separable and $||| \cdot ||||$ is equivalent to the original norm
on $Y$). For $m=1,2,...$ put $g_m = f_{n_{k_{i_m}}}$. Since $\{
g_m - f_0\}$ is a $w^*-$basic sequence, it follows that for some
sequence $\{ x_m \}_{m=1}^\infty$ in $X$

\begin{equation}\label{bo}
\{ (g_m - f_0,x_m) \}_{m=1}^\infty \ \ \ \ {\rm is \ a \
biorthogonal \ sequence.}
\end{equation}

\noindent Only two cases must be now considered.

1) For some sequence $\{ m_j\}_{j=1}^\infty$ we have $f_0
(x_{m_j})=0, \, j=1,2,...$: in this case $\{(g_{m_j},x_{m_j})\}$
would be a biorthogonal system, contradicting the fact that
$\{g_{m_j} \}$ is almost overtotal on $X$.

2) There exists $q$ such that for any $m\geq q$ we have $f_0
(x_m)\not = 0.$ For any $j >q$, from (\ref{bo}) it follows
$$0=(g_{3j}-f_0)(f_0(x_{3j-1}) \cdot x_{3j-2}-f_0(x_{3j-2}) \cdot x_{3j-1})
=g_{3j}(f_0(x_{3j-1}) \cdot x_{3j-2}-f_0(x_{3j-2}) \cdot
x_{3j-1}).$$

\noindent It follows that the almost overtotal sequence
$\{g_{3j}\}_{j=q}^{\infty}$ annihilates the subspace $W =[\{
f_0(x_{3j-1}) \cdot x_{3j-2}-f_0(x_{3j-2}) \cdot
x_{3j-1}\}_{j=q}^\infty] \subset X$: being $\{ x_m
\}_{m=1}^\infty$ a linearly independent sequence, $W$ is
infinite-dimensional, a contradiction.

Hence (\ref{not to}) does not work and we are done. \ $\black$

\bigskip

As an immediate consequence we get the following Corollary.

\begin{corollary}\label{coro}
Let $X$ be an infinite-dimensional Banach space and
$\{f_n\}\subset B_{X^*}$ be a sequence that is not relatively
norm-compact. Then there is an infinite-dimensional subspace
$Y\subset X$ such that $|\{f_n\}\cap Y^{\perp}|=\infty .$ For
instance this is true for any $\delta$-separated sequence
$\{f_n\}\subset B_{X^*} \ (\delta
>0)$.
\end{corollary}

\bigskip

\section{Applications}

The following theorem easily follows from Corollary \ref{coro}.
\begin{theorem}\label{anni}
Let $X\subset C(K)$ be an infinite-dimensional subspace of $C(K)$
where $K$ is metric compact. Assume that, for $\{t_n\}_{n \in
\mathbb N} \subset K$, the sequence $\{t_n |_X \} \subset X^*$ is
not relatively norm-compact. Then there are an
infinite-dimensional subspace $Y\subset X$ and a subsequence
$\{t_{n_k}\}_{k \in \mathbb N}$ such that $y(t_{n_k})=0$ for any
$y\in Y$ and for any $k \in \mathbb N$.
\end{theorem}

\noindent {\bf Remark.} \ Sequences $\{t_n\}\subset K$ as required
in the statement of Theorem \ref{anni} always exist: trivially,
for any sequence $\{t_n\}$ dense in $K$, the sequence $\{t_n
|_X\}$, being a 1-norming sequence for $X$, cannot be relatively
norm-compact (since $X$ is infinite-dimensional).

In \cite{EGS} the Authors proved that, for any
infinite-dimensional subspace $X\subset C(K)$, there are an
infinite-dimensional subspace $Y\subset X$ and a sequence
$\{t_k\}_{k \in \mathbb N} \subset K$ such that $y(t_k)=0$ for any
$y\in Y$ and any $k \in \mathbb N$. Theorem \ref{anni} strengthens
this result. In fact actually, for any infinite-dimensional
subspace $X\subset C(K)$, we can find such a sequence $\{ t_k\}$
as a suitable subsequence $\{ w_{n_k}\}$ of any prescribed
sequence $\{ w_n \} \subset K$ for which $\{ w_n |_X \} \subset
X^*$ is not relatively norm-compact.

\bigskip

In 2003 R. Aron and V. Gurariy asked: does there exist an
infinite-dimensional subspace of $l_{\infty}$ every non-zero
element of which has only finitely many zero-coordinates? Let us
reformulate this question in the following equivalent way: {\it
does there exist an infinite-dimensional subspace $Y \subset
l_\infty$ such that the sequence $\{ e_n|_Y \}$ of the
``coordinate functionals'' is overtotal on $Y$?}

Since the sequence $\{ e_n|_Y \}$ is norming for $Y$, it is not
norm-compact ($Y$ is infinite-dimensional), hence by Theorem
\ref{compots} it cannot be overtotal on $Y$. So the answer to the
Aron-Gurariy's question is negative. Actually we can say much
more. In fact, from Theorem \ref{compots} it follows that there
exist an infinite-dimensional subspace $Z \subset Y$ and a
strictly increasing sequence $\{n_k \}$ of integers such that $\{
e_{n_k}(z) = 0 \}$ for every $z \in Z$ and $k \in \mathbb N$.

Note that the Aron-Gurariy's question was answered via a different
argument in \cite{CS}.

\bigskip

The next Theorem generalizes the previous argument.

\begin{theorem}
\label{embedding} Let $X$ be a separable infinite-dimensional
Banach space and $T: X \to l_\infty$ be a one-to-one bounded non
compact linear operator. Then there exist an infinite-dimensional
subspace $Y\subset X$ and a strictly increasing sequence $\{n_k\}$
of integers  such that $e_{n_k}(Ty)=0$ for any $y\in Y$ and for
any $k$ ($e_n$ the ``$n$-coordinate functional'' on $l_\infty$).
\end{theorem}

\noindent{\bf Proof.} Assume to the contrary that for any sequence
of integers $\{n_k\}$ we have

\noindent $\dim\bigl( \{T^* e_{n_k}\}^{\top} \bigr) < \infty .$
Then the sequence $\{T^*(e_n)\} \subset X^*$ is almost overtotal
on $X$, so $K = || \cdot ||-{\rm cl}\{T^*e_n\}$ is norm-compact in
$X^*$ by Theorem \ref{compots}. Clearly we can consider $B_X$ as a
subset of ${\cal C}(K)$ (by putting, for $x \in B_X$ and $t \in
K$, $x(t) = t(x)$). We claim that $B_X$ is relatively norm-compact
in ${\cal C}(K)$. In fact, $B_X$ is clearly bounded in ${\cal
C}(K)$ and its elements are equi-continuous since, for $t_1, t_2
\in K$ and $x \in B_X$, we have
$$|x(t_1) - x(t_2)| \leq ||x|| \cdot ||t_1 - t_2|| \leq ||t_1 - t_2||:$$
we are done by the Ascoli-Arzel\`a theorem. Since, for $x \in X$
we have $||x||_{{\cal C}(K)} = ||Tx||_{l_\infty}$, $T(B_X)$ is
relatively norm-compact in $l_\infty$ too. This leads to a
contradiction since we assumed that $T$ is not a compact operator.
\ $\black$

\bigskip

Let now $X$ be an infinite-dimensional space and $\{ f_n\} \subset
X^*$ a norming sequence for $X$. By Theorem \ref{compots}, the
fact that $\{ f_n\}$ is not relatively norm-compact immediately
forces $\{ f_n\}$ not to be overtotal on $X$. Since any norming
sequence is a total sequence, it follows that any norming sequence
for any infinite-dimensional space $X$ admits some subsequence
that is not a norming sequence for $X$. In other words and
following our terminology, ``overnorming'' sequences do not exist.

\bigskip

As one more application of Theorem \ref{compots} we obtain the
following Theorem.

\begin{theorem}\label{ttt}
\label{l_p} Let $X, \ Y$ be infinite-dimensional Banach spaces,
$Y$ having an unconditional basis $\{u_i \}_{i=1}^\infty$ with
$\{e_i \}_{i=1}^\infty$ as the sequence of the associated
coordinate functionals. Let $T: X \to Y$ be a one-to-one bounded
non compact linear operator. Then there exist an
infinite-dimensional subspace $Z\subset X$ and a strictly
increasing sequence $\{k_l\}$ of integers such that
$e_{k_l}(Tz)=0$ for any $z\in Z$ and any $l \in \mathbb N$.
\end{theorem}

To prove Theorem \ref{l_p} we need some preparation. First note
that, without loss of generality, from now on we may assume that
$T$ has norm one and that the unconditional basis $\{u_i
\}_{i=1}^\infty$ is normalized and unconditionally monotone (i.e.,
if $x = \sum_{i=1}^\infty \alpha_iu_i$ and $\sigma \subset \mathbb
N$, then $||\sum_{i \in \sigma}\alpha_i\beta_iu_i|| \leq ||x||$
for any choice of $\beta_i$ with $|\beta_i| \leq 1$; see for
instance \cite{Si1}, Theorem 17.1).

\medskip

In the proof of Theorem \ref{l_p} we will use the following two
technical Lemmas.

\begin{lemma}
\label{lemma 1 l_p} Let $X$, $Y$ and $T$ be as in the statement of
Theorem \ref{l_p}. Then there exists $\delta > 0$ such that, for
any natural integer $m$, some point $z \in B_X$ exists (depending
on $m$) such that $||Tz|| \geq \delta$ and the first $m$
coordinates of $Tz$ are 0.
\end{lemma}

\noindent{\bf Proof.} Let us start by proving that, keeping
notation as in the statement of Theorem \ref{l_p},
\begin{equation}\label{fact}
\exists \{x_k \}_{k=1}^\infty \subset B_X, \ \exists \ 0 < \beta <
1: e_i(Tx_k) \to 0\ {\rm as}\ k \to \infty \ \forall i \in \mathbb
N\ \wedge \ ||Tx_k||
> \beta \ \forall k \in \mathbb N.
\end{equation}

In fact, let $\{z_n \}_{n=1}^\infty$ be any $r-$separated sequence
in $T(B_X)$ for some $r >0$ ($T(B_X)$ is not pre-compact). By a
standard diagonal procedure we can select a subsequence $\{
z_{n_k}\}$ such that, for any $i \in \mathbb N$, the numbers
$e_i(z_{n_k})$ converge as $k \to \infty$. Of course, for any $i$
we have $e_i(z_{n_k} - z_{n_{k+1}}) \to 0$ as $k \to \infty$ with
$||z_{n_k} - z_{n_{k+1}}|| \geq r$. For each $k$, put $2y_k =
z_{n_{2k}} - z_{n_{2k+1}}$: since $T(B_X)$ is both convex and
symmetric with respect to the origin, it is clear that $\{
y_k\}_{k=1}^\infty \subset T(B_X)$ too; moreover for any $k$ we
have $||y_k|| > r/2$ and for any $i$ we have $e_i(y_k) \to 0$ as
$k \to \infty$. So it is enough to assume $x_k = T^{-1}y_k$ for
any $k$ and $\beta = r/2$ and (\ref{fact}) is proved.

Now fix $m \in \mathbb N$. Put $L = {[\{ T^*e_n\}_{n=1}^m}]^\top$
and let $x \in X$. Then, denoting by $q: X \to X/L$ the quotient
map, for some positive constant $C_m$ independent on $x$ it is
true that
$${\rm dist}(x,L) = ||q(x)|| = {\rm Sup}\{ |f(q(x))|: f \in S_{(X/L)^*} \} =
$$
\begin{equation}\label{quotient}
= {\rm Sup}\{ |g(x)|: g \in S_{[\{ T^*e_n\}_{n=1}^m]} \} \leq
C_m{\rm Max}\{ |e_n(Tx)|: 1 \leq n \leq m \}.
\end{equation}

\noindent Take $\{x_k \}_{k=1}^\infty$ as in (\ref{fact}): some
$\tilde k \in \mathbb N$ exists such that
$$ C_m{\rm Max}\{ |e_n(Tx_{\tilde k})|: 1 \leq n \leq m \} < \beta /2$$
that by (\ref{quotient}) implies
$${\rm dist}(x_{\tilde k},L) < \beta /2.$$
Let $2z \in L$ be such that $||x_{\tilde k} - 2z|| < \beta /2$:
clearly $||z|| < 1$ and $\displaystyle ||Tz|| > (||Tx_{\tilde k}||
- \beta /2)/2$, so, since $||Tx_{\tilde k}|| > \beta$, we are done
by assuming $\delta = \beta /4$. \ $\black$

\bigskip

\begin{lemma}
\label{lemma 2 l_p} Let $Y$ be as in the statement of Theorem
\ref{l_p}. Then for any
\begin{equation}\label{vw}
n \in \mathbb N, \ \ 0 < \varepsilon \leq 1/2, \ \ v =
\sum_{i=1}^n v_iu_i, \ \ w = \sum_{i=1}^n w_iu_i \ \ with \ \
||v|| < \varepsilon^2 \ \ and \ \ ||w|| > 1 - \varepsilon,
\end{equation}
there exists $j, \ 1 \leq j \leq n$, such that $|v_j| <
\varepsilon|w_j|$.
\end{lemma}

\noindent{\bf Proof.} Recall that, under our assumptions, basis
$\{ u_i\}$ is unconditionally monotone. Hence, without loss of
generality, we may assume that $w_i\not = 0, i=1,...,n$. Moreover,
for any $n \in \mathbb N$, any scalars $\alpha_1,...,\alpha_n$ and
$|\beta_1|,...,|\beta_n| \leq 1$, the following is true
\begin{equation}\label{constant}
||\sum_{i=1}^n \beta_i\alpha_iu_i|| \leq ||\sum_{i=1}^n
\alpha_iu_i||.
\end{equation}
Assume to the contrary that some $v,w$ exist satisfying (\ref{vw})
for some $\varepsilon$, $0 < \varepsilon < 1/2$, for which $|v_i|
\geq \varepsilon |w_i|$ (i.e. $\varepsilon |w_i/v_i| \leq 1$) for
every $i$, $1 \leq i \leq n$. By putting in (\ref{constant})
$\alpha_i = v_i/\varepsilon$ and $\beta_i = \varepsilon w_i/v_i$
for any $i$, we get
$$1 - \varepsilon < ||\sum_{i=1}^n w_iu_i|| \leq
||\sum_{i=1}^n v_iu_i/\varepsilon || <  \epsilon$$ that gives
$\varepsilon > 1/2$, a contradiction. \ $\black$

\bigskip

\noindent{\bf Proof of Theorem  \ref{l_p}.}

By Lemma \ref{lemma 1 l_p}, a bounded sequence $\{x_n
\}_{n=1}^\infty, \ ||x_n|| < R$ for some $R > 0$, can be found in
$X$ such that $Tx_n \in S_Y$ for every $n$ and $e_j(Tx_n) = 0, \
j=1,...,n$. For any $n$ put
$$Tx_n = y_n = \sum_{i=n+1}^\infty y_n^i u_i.$$

Now we are going to construct a subsequence $\{y_{n_k}
\}_{k=0}^\infty$ of $\{y_n \}_{n=1}^\infty$ with special
properties.

Put in short $1/2^{n+1} = \varepsilon_n, \ n=1,2,..\, .$.

Put $n_0 = 1$ and let $p_0 > n_0$ be such that $y_{n_0}^{p_0} \neq
0$.

Take $n_1 \geq p_0$ such that
$$||\sum_{n_1+1}^\infty y_{n_0}^iu_i|| < \varepsilon_1^2.$$
Let $n_2 > n_1$ such that (remember that our basis is
unconditionally monotone)
$$||\sum_{n_2+1}^\infty |y_{n_0}^i|u_i||
+ ||\sum_{n_2+1}^\infty |y_{n_1}^i|u_i|| < \varepsilon_2^2$$ and
consider the two vectors
$$v_1 = \sum_{n_1+1}^{n_2} y_{n_0}^iu_i , \ \ \ \ \ \
w_1 = \sum_{n_1+1}^{n_2} y_{n_1}^iu_i:$$ clearly we have $||v_1||
< \varepsilon_1^2$ and $||w_1|| > 1 - \varepsilon_2^2 > 1 -
\varepsilon_1$, hence by Lemma \ref{lemma 2 l_p} an integer $p_1,
\ n_1+1 \leq p_1 \leq n_2$, can be found such that
$${{|y_{n_0}^{p_1}|} \over {|y_{n_1}^{p_1}|}} < \varepsilon_1.$$
Now take $n_3 > n_2$ such that
$$\sum_{j=0}^2\sum_{n_3+1}^\infty ||\, |y_{n_j}^i|u_i|| < \varepsilon_3^2$$
and consider the two vectors
$$v_2 = \sum_{n_2+1}^{n_3} (|y_{n_0}^i| + |y_{n_1}^i|)u_i , \ \ \ \ \ \
w_2 = \sum_{n_2+1}^{n_3} y_{n_2}^iu_i:$$ clearly we have $||v_2||
< \varepsilon_2^2$ and $||w_2|| > 1 - \varepsilon_3^2 > 1 -
\varepsilon_2$, hence by Lemma \ref{lemma 2 l_p} an integer $p_2,
\ n_2+1 \leq p_1 \leq n_3$, can be found such that
$${{|y_{n_0}^{p_2}| + |y_{n_1}^{p_2}|} \over {|y_{n_2}^{p_2}|}} < \varepsilon_2.$$
It is now clear how to iterate the process, so getting a sequence
$\{y_{n_k}\}_{k=0}^\infty$ in $S_{T(X)}$, a corresponding
subsequence $\{ p_k \}_{k=0}^\infty$ being determined such that
for $k \geq 0$
\begin{equation}\label{n_k}
n_k+1 \leq p_k \leq n_{k+1} \ \ \ \wedge \ \ \
{{\sum_{j=0}^{k-1}|y_{n_j}^{p_k}|} \over {|y_{n_k}^{p_k}|}}  <
\varepsilon_k.
\end{equation}

Put
$$E = [\{ y_{n_k}\}_{k=0}^\infty], \ \ \ W = T^{-1}(E), \ \ \
\tilde e_k = e_{p_k}/y_{n_k}^{p_k},\ k=0,1,2,...\ .$$

Clearly we have
\begin{equation}\label{tilde e}
\tilde e_k(y_{n_i}) = 0 \ \ \ {\rm if} \ k < i, \ \ \ \ \ \ \tilde
e_k(y_{n_k}) = 1, \ \ \ \ \ k = 0,1,2,... \ .
\end{equation}

Note that, by our construction, $\{y_{n_k}\}_{k=0}^\infty$ is a
sufficiently small perturbation of a block basis of the basis
$\{u_i\}.$ Hence it is an unconditional basis for $E$. Let $B$ its
basis constant.

We claim that $\{ T^*\tilde e_k|_W \}_{k=1}^\infty \subset W^*$ is
a bounded sequence. Clearly it is enough to prove that $\{ \tilde
e_k|_E \}_{k=1}^\infty$ is bounded. In fact, for any $k \in
\mathbb N$ and any $y = \sum_{i=0}^\infty a_i y_{n_i} \in S_E$,
taking into account (\ref{tilde e}) and (\ref{n_k}) we have
$$|\tilde e_k(y)| = |\tilde e_k(\sum_{i=0}^\infty a_iy_{n_i})| =
|\tilde e_k(\sum_{i=0}^k a_iy_{n_i})| \leq \sum_{i=0}^k
|a_i||\tilde e_k(y_{n_i})| \leq 2B(\varepsilon_k +1) < 4B.$$

Moreover we claim that it is a $1/R$-separated sequence. In fact
for any $k,m$ with $k > m \geq 0$, again remembering (\ref{tilde
e}), we have
$$||T^*\tilde e_k - T^*\tilde e_m|| \geq |(T^*\tilde e_k)(x_{n_k}/R) -
(T^*\tilde e_m)(x_{n_K}/R)| =$$
$$= (1/R)|(\tilde e_k)(y_{n_k}) -
(\tilde e_m)(y_{n_k})| = 1/R.$$

Hence, by Theorem \ref{compots}, the sequence $\{ T^*\tilde e_k|_W
\}_{k=1}^\infty$ cannot be almost overtotal on $W$: it means that
there is an infinite-dimensional subspace $Z \subset W$ that
annihilates some subsequence of the sequence $\{ T^*\tilde e_k
\}$.

The proof is complete. \ $\black$

\bigskip

\noindent{\bf Remark.} \ The particular case of Theorem \ref{ttt}
when $X \subset l_p, \ Y = l_p \ (1 \leq p < \infty), \ T = {\rm
Id}_{|X}$ was proved in \cite{CS}.

\bigskip

Vladimir P. Fonf

Department of Mathematics

Ben-Gurion University of the Negev

84105 Beer-Sheva, Israel

E-mail address: fonf@math.bgu.ac.il

\bigskip
\bigskip

Clemente Zanco

Dipartimento di Matematica

Universit\`a degli Studi

Via C. Saldini, 50

20133 Milano MI, Italy

E-mail address: clemente.zanco@unimi.it

ph. ++39 02 503 16164 \ \ \ \ fax  ++39 02 503 16090

\end{document}